# Studies on Frequency Response Optimized Integrators Considering Second Order Derivative

Sheng Lei, *Student Member, IEEE*, and Alexander Flueck, *Senior Member, IEEE*

*Abstract*—This paper presents comprehensive studies on frequency response optimized integrators considering second order derivative regarding their numerical error, numerical stability and transient performance. Frequency domain error analysis is conducted on these numerical integrators to reveal their accuracy. Numerical stability of the numerical integrators is investigated. Interesting new types of numerical stability are recognized. Transient performance of the numerical integrators is defined to qualitatively characterize their ability to track fast decaying transients. This property is related to unsatisfactory phenomena such as numerical oscillation which frequently appear in time domain simulation of circuits and systems. Transient performance analysis of the numerical integrators is provided. Theoretical observations from the analysis of the numerical integrators are verified via time domain case studies.

## I. INTRODUCTION

TIME domain simulation forms a powerful tool for computer-aided analysis and design of electronic circuits [1]-[4] and electric power systems [5]-[6]. When performing time domain simulation, one always prefers accuracy and efficiency [7]-[9]. Accuracy can be understood as the degree to which the computed results match reality or the analytical results. Efficiency is inversely proportional to the total time consumption required for a simulation run. For a given simulation method or scheme, in order to achieve higher accuracy, the required step size typically has to be shrunk [7], [10]. As a result, more time steps are required for a simulation run, leading to higher computational burden and increased time consumption.

In fact, both objectives may be achieved simultaneously as the required, or an even better, degree of accuracy may be obtained with less time consumption by highly accurate simulation methods or schemes, which enable larger step sizes, resulting in fewer time steps and less computation [7], [11]-[15]. Specifically, in time domain simulation of electronic circuits and electric power systems, the error in discretizing the ordinary differential equations (ODEs) describing the dynamics of the system under study mainly undermines the overall accuracy of the simulation [7]-[9], [14]-[15]. Therefore much research effort has been directed to introducing highly accurate numerical integrators to the discretization process. In [7], [11]-[13], [16]-[17], the Obreshkov (also spelled as Obrechkoff) family of numerical integrators, considering higher order derivatives of the differential state variables, are introduced into time domain simulation of electronic circuits.

The accuracy in the discretization process may be improved from another angle by taking into account the dominant frequency component of the signals. In certain situations, signals inside the studied system are dominated by a specific frequency. For electronic communication circuits, information as slow variants are modulated onto the high-frequency carrier signals so that band-pass signals are created [8]-[9], [18]. The frequency spectrum of these signals is concentrated around the carrier frequency within a narrow bandwidth. For power system stability studies, voltage and current waveforms in the network are dominated by the nominal fundamental frequency (50 or 60 Hz), after fast transients have died down which are not of interest [9], [19]-[20]. Note that the state variable of some ODEs may belong to these signals. If the numerical integrators are designed so that they introduce only negligible error to the dominant frequency component, namely if they are frequency response optimized, the overall accuracy of the simulation will be improved. This frequency response optimization idea of numerical integrators is proposed in [14] and applied to modify the well-known implicit trapezoidal method, where the frequency response of the error is set to zero at the nominal fundamental frequency. The modified implicit trapezoidal method is then adopted in power system transient simulation [14]-[15].

Extending the frequency response optimization idea in [14]-[15] to the Obreshkov-like numerical integrators, [21] proposes frequency response optimized integrators considering second order derivative. Based on the proposed numerical integrators and some others from the original Obreshkov family, [21] further puts forward a novel power system transient simulation scheme. Note that external inputs to power systems may not have explicit analytical expressions. For example, wind power and sun radiation as energy sources are well known to be stochastic [22]-[23] and are often given as time series for

The authors are with the Department of Electrical and Computer Engineering, Illinois Institute of Technology, Chicago, IL 60616 USA (e-mail: slei3@hawk.iit.edu; flueck@iit.edu).

Sheng Lei is also with the Mathematics and Computer Science Division, Argonne National Laboratory, Lemont, IL 60439 USA.

simulation [24]-[26]. As a result, analytical calculation of their derivative and higher order derivatives is impossible. Numerical calculation of this information [27]-[28], however, introduces non-negligible error, which will deteriorate the overall accuracy, even though Obreshkov-like numerical integrators are used. In that case, considering higher order derivatives becomes an unnecessary burden without improving the accuracy. With this constraint in mind, derivatives higher than the second order are not considered in [21]. The simulation scheme is shown to considerably enhance the computational efficiency while providing the desired accuracy [21]. Unfortunately, properties of the underlying numerical integrators are not studied in [21] due to space limitation.

In this paper, comprehensive studies are given on frequency response optimized integrators considering second order derivative. Note that some numerical integrators from the original Obreshkov family are also understood as this type in the sense that they are optimized for signals around 0 Hz (slow variants) and they do consider second order derivative [21]. Contributions of this paper include:
1) Frequency response of the error of the numerical integrators is given. The conventional Taylor expansion-based truncation error analysis [27]-[30] is not able to reveal the accuracy of numerical integrators for frequency components other than 0 Hz. Therefore error analysis is performed in the frequency domain in this paper following [14]-[15], [31].
2) Numerical stability of the numerical integrators is analyzed. Interesting new types of numerical stability are recognized.
3) Transient performance of the numerical integrators is defined to qualitatively characterize their accuracy for fast and monotonically decreasing transients. This property is related to unsatisfactory phenomena in transient simulation such as numerical oscillation [32]-[33] induced by numerical integrators. Transient performance of the numerical integrators is analyzed.
4) Error analysis and transient performance of the numerical integrators are verified by case studies from the time domain perspective.

The rest of the paper is organized as follows. Section II reviews formulas of frequency response optimized integrators considering second order derivative and gives frequency domain error analysis. Section III analyzes numerical stability of the numerical integrators. Transient performance of the numerical integrators is analyzed in Section IV. Time domain case studies verify error analysis and transient performance of the numerical integrators in Section V. Finally, Section VI concludes the paper and points out some directions for future research.

## II. FREQUENCY RESPONSE OPTIMIZED INTEGRATORS CONSIDERING SECOND ORDER DERIVATIVE

Consider a general ODE of the form

$$\dot{x} = f(t, x, u) \quad (1)$$

where $t$ denotes the time instant; $x$ denotes the state variable; $u$ denotes the input; $f$ is a function depending on $t$, $x$ and $u$. A single input is considered here for clear presentation, but the same idea applies to multiple inputs.

Applying a numerical integrator considering second order derivative to discretize (1), the following algebraic equation is obtained at $t$ [7], [21], [29]-[30]

$$x_t = a_{-1} x_{t-h} + b_0 \dot{x}_t + b_{-1} \dot{x}_{t-h} + c_0 \ddot{x}_t + c_{-1} \ddot{x}_{t-h} \quad (2)$$

where $h$ is the step size; $a_{-1}$, $b_0$, $b_{-1}$, $c_0$ and $c_{-1}$ are coefficients to be determined. A specific selection for these coefficients determines a numerical integrator.

The required second order derivative of $x$ in (2) can be calculated by taking derivative on both sides of (1)

$$\ddot{x} = \frac{\partial f}{\partial t} + \frac{\partial f}{\partial x} \dot{x} + \frac{\partial f}{\partial u} \dot{u} = \frac{\partial f}{\partial t} + \frac{\partial f}{\partial x} f(t,x,u) + \frac{\partial f}{\partial u} \dot{u} \quad (3)$$

Obviously the derivative of $u$ is required. If the expression of $u$ is given, this information can be calculated analytically. However, in the general setting, the expression is unavailable and the derivative has to be calculated by numerical derivative [27]-[28].

Reference [14] proposes modifying the implicit trapezoidal method so that the frequency response of the error is set to zero at a specified nonzero angular frequency $\omega_{select}$. Such a modification makes the resulting numerical integrator especially accurate for those ODEs of which the state variable is dominated by the frequency component around $\omega_{select}$. This idea can be extended to numerical integrators considering higher order derivatives. In particular, performing the Laplace transform on both sides of (2)

$$X = a_{-1} X e^{-sh} + b_0 s X + b_{-1} s X e^{-sh} + c_0 s^2 X + c_{-1} s^2 X e^{-sh} \quad (4)$$

The $s$-domain error expression is

$$X - (a_{-1} X e^{-sh} + b_0 s X + b_{-1} s X e^{-sh} + c_0 s^2 X + c_{-1} s^2 X e^{-sh}) \quad (5)$$

The $s$-domain relative error expression is

$$1 - (a_{-1} e^{-sh} + b_0 s + b_{-1} s e^{-sh} + c_0 s^2 + c_{-1} s^2 e^{-sh}) \quad (6)$$

If the coefficients are chosen so that 0 is a root of (6), then the frequency response of the error is set to zero at 0 Hz, leading to a numerical integrator which is accurate for slow variants. Similarly, if $j\omega_{select}$ and $-j\omega_{select}$ are made a pair of roots of (6) by properly choosing the coefficients, then the frequency response of the error is set to zero at the angular frequency $\omega_{select}$, making the numerical integrator accurate for signals with a dominant frequency component at $\omega_{select}$. A multiple root introduces smaller error around the specified frequency. Desirable selection is achieved by solving equations regarding these root conditions.

### A. Making $j\omega_{select}$ and $-j\omega_{select}$ a Single Root and 0 a Triple Root

Reference [21] proposes a numerical integrator expressed as

$$a_{-1}=1, \quad b_0=\frac{h}{2}, \quad b_{-1}=\frac{h}{2}$$
$$c_0 = -\frac{1}{\omega_{select}^2} + \frac{h}{2\omega_{select}}\cot(\frac{\omega_{select}h}{2}) \quad (7)$$
$$c_{-1} = \frac{1}{\omega_{select}^2} - \frac{h}{2\omega_{select}}\cot(\frac{\omega_{select}h}{2})$$

With this set of coefficients, $j\omega_{select}$ and $-j\omega_{select}$ are a single root and 0 is a triple root of the relative error expression (6). This numerical integrator is referred to as Integrator A hereafter. Integrator A is accurate for the frequency component around $\omega_{select}$; it is also rather accurate for slow variants. Magnitude frequency response of the error of Integrator A is plotted with different step sizes in Fig. 1, where $\omega_{select}$ is set at 60 Hz.

### B. Making $j\omega_{select}$, $-j\omega_{select}$ and 0 a Single Root Respectively while Letting $b_{-1}$ and $c_{-1}$ Be 0

Reference [21] proposes a numerical integrator expressed as

$$a_{-1}=1, \quad b_0=\frac{\sin(\omega_{select}h)}{\omega_{select}}, \quad b_{-1}=0 \quad (8)$$
$$c_0=\frac{\cos(\omega_{select}h)-1}{\omega_{select}^2}, \quad c_{-1}=0$$

With this set of coefficients, $j\omega_{select}$, $-j\omega_{select}$ and 0 are a single root of the relative error expression (6) respectively. As $b_{-1}$ and $c_{-1}$ are both 0, it is suitable for dealing with discontinuities. This numerical integrator is referred to as Integrator B hereafter. Integrator B is accurate for the frequency component around $\omega_{select}$ and 0 Hz. Magnitude frequency response of the error of Integrator B is plotted with different step sizes in Fig. 2, where $\omega_{select}$ is set at 60 Hz.

### C. Other Numerical Integrators Considering Second Order Derivative

A numerical integrator considering second order derivative from the Obreshkov family is [29]-[30], [34]

$$a_{-1}=1, \quad b_0=\frac{h}{2}, \quad b_{-1}=\frac{h}{2}, \quad c_0=-\frac{h^2}{12}, \quad c_{-1}=\frac{h^2}{12} \quad (9)$$

With this set of coefficients, 0 is a quintuple root of the relative error expression (6), implying that this numerical integrator is highly accurate for slow variants. This numerical integrator is referred to as Integrator C hereafter. Magnitude frequency response of the error of Integrator C is plotted with different step sizes in Fig. 3.

An implicit second order Taylor series method [27], [29] is also from the Obreshkov family, which is express as

$$a_{-1}=1, \quad b_0=h, \quad b_{-1}=0, \quad c_0=-\frac{h^2}{2}, \quad c_{-1}=0 \quad (10)$$

With this set of coefficients, 0 is a triple root of the relative error expression (6), implying that this numerical integrator is rather accurate for slow variants. As $b_{-1}$ and $c_{-1}$ are both 0, it is suitable for dealing with discontinuities. This numerical integrator is referred to as Integrator D hereafter. Magnitude frequency response of the error of Integrator D is plotted with different step sizes in Fig. 4.

### D. Remarks

According to the frequency response of the error, if the state variable of an ODE is a pure sinusoid at $\omega_{select}$, then Integrators A and B are exactly accurate for it despite the step size, as long as the coefficients are not made infinite.

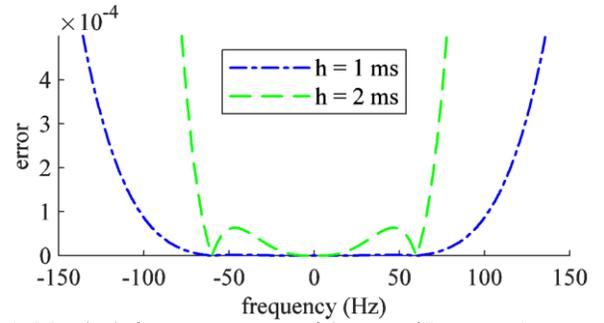

Fig. 1. Magnitude frequency response of the error of Integrator A.

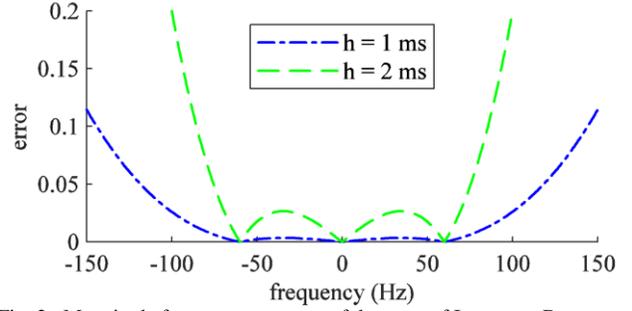

Fig. 2. Magnitude frequency response of the error of Integrator B.

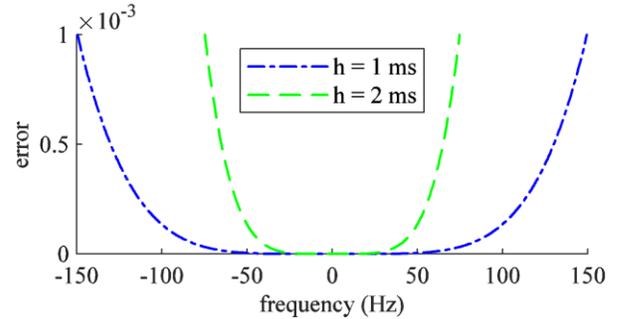

Fig. 3. Magnitude frequency response of the error of Integrator C.

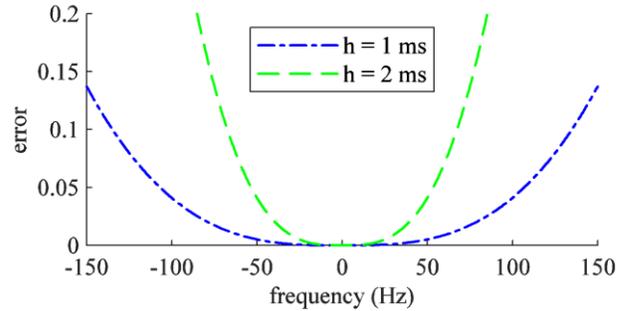

Fig. 4. Magnitude frequency response of the error of Integrator D.

From the frequency response of the error of Integrators C and D, it is observed that they introduce considerable error at the nominal fundamental frequency of power systems if a step size of 1 or 2 ms is adopted. In order to reduce the error, a smaller step size has to be used than that which is sufficient for depicting the waveform, making the computation inefficient. Instead, Integrators A and B should be used.

For Integrators C and D, the multiplicity of the root at 0 of the relative error expression (6) coincides with the order of their local truncation error [29]-[30], [34]. The relation between the conventional Taylor expansion-based truncation error analysis and the frequency domain error analysis may be an interesting direction for further research.

## III. NUMERICAL STABILITY

Consider the following test ODE with the given initial value condition

$$\begin{cases} \dot{x} = \lambda x \\ x(0) = 1 \end{cases} \quad (11)$$

where $\lambda$ is a complex-valued number with negative real part, which is understood as the eigenvalue of the equation (11); at the time instant 0, the value of $x$ is 1. Clearly $x$ tends to 0 as $t$ tends to infinity [35]-[36].

Now a numerical integrator is applied to solve (11) given a step size $h$ which is real-valued and positive. The ODE is discretized by the numerical integrator to be converted into an algebraic equation at each time step. It is expected that the solution to the series of algebraic equations, namely the discretized ODE, also tends to 0 as the time step index tends to infinity. Otherwise the numerical solution does not match the trend of the analytical solution and is considered inaccurate. It is shown that for the numerical solution from numerical integrators such as the forward Euler method to match the trend, the step size has to be very small if the magnitude of $\lambda$ is large [27]-[30], which causes computational inefficiency; otherwise the numerical solution will blow up. The eigenvalues of a system under study are generally unknown in advance, making the step size selection difficult if this type of numerical integrator is adopted. It is desirable for the numerical solution to match the trend despite the step size so that the constraint on step size selection from this aspect can be relaxed. Such a satisfactory property of numerical integrators is called A-stability [27]-[30]. Note that $\lambda$ with large magnitude corresponds to fast transients or high frequency oscillations [35]. Therefore A-stability qualitatively characterizes the accuracy of numerical integrators for these types of signals.

Integrators C and D have already been proven to be A-stable in the literature [7], [27], [30], [34]. Moreover Integrator D is in fact L-stable [7], [27], [30]. Therefore this paper does not need to discuss their numerical stability. Numerical stability of Integrators A and B is discussed in this section.

### A. Numerical Stability of Integrator A

Applying Integrator A, the test ODE (11) is discretized as

$$x_t = x_{t-h} + b_0 \dot{x}_t + b_{-1} \dot{x}_{t-h} + c_0 \ddot{x}_t + c_{-1} \ddot{x}_{t-h} \\ = x_{t-h} + b_0 \lambda x_t + b_{-1} \lambda x_{t-h} + c_0 \lambda^2 x_t + c_{-1} \lambda^2 x_{t-h} \quad (12)$$

For the numerical solution to tend to 0, it is required that

$$\left| \frac{1 + b_{-1}\lambda + c_{-1}\lambda^2}{1 - b_0\lambda - c_0\lambda^2} \right| < 1 \quad (13)$$

which is guaranteed by

$$|1 + b_{-1}\lambda + c_{-1}\lambda^2| < |1 - b_0\lambda - c_0\lambda^2| \quad (14)$$

guaranteed by

$$|1 + b_{-1}\lambda + c_{-1}\lambda^2|^2 < |1 - b_0\lambda - c_0\lambda^2|^2 \quad (15)$$

or

$$(1 + b_{-1}\lambda + c_{-1}\lambda^2)(1 + b_{-1}\lambda + c_{-1}\lambda^2)^* \\ < (1 - b_0\lambda - c_0\lambda^2)(1 - b_0\lambda - c_0\lambda^2)^* \quad (16)$$

where * denotes the complex conjugate. (16) is equivalent to

$$(1 + b_{-1}\lambda + c_{-1}\lambda^2)(1 + b_{-1}\lambda^* + c_{-1}\lambda^{*2}) \\ < (1 - b_0\lambda - c_0\lambda^2)(1 - b_0\lambda^* - c_0\lambda^{*2}) \quad (17)$$

namely

$$\begin{array}{ll} (1 + b_{-1}\lambda^* + c_{-1}\lambda^{*2}) & (1 - b_0\lambda^* - c_0\lambda^{*2}) \\ + b_{-1}\lambda(1 + b_{-1}\lambda^* + c_{-1}\lambda^{*2}) & < -b_0\lambda(1 - b_0\lambda^* - c_0\lambda^{*2}) \\ + c_{-1}\lambda^2(1 + b_{-1}\lambda^* + c_{-1}\lambda^{*2}) & -c_0\lambda^2(1 - b_0\lambda^* - c_0\lambda^{*2}) \end{array} \quad (18)$$

namely

$$\begin{array}{l} 1 + b_{-1}\lambda^* + c_{-1}\lambda^{*2} \\ + b_{-1}\lambda + b_{-1}^2|\lambda|^2 + b_{-1}c_{-1}|\lambda|^2\lambda^* \\ + c_{-1}\lambda^2 + c_{-1}b_{-1}|\lambda|^2\lambda + c_{-1}^2|\lambda|^4 \\ < \\ 1 - b_0\lambda^* - c_0\lambda^{*2} \\ - b_0\lambda + b_0^2|\lambda|^2 + b_0c_0|\lambda|^2\lambda^* \\ - c_0\lambda^2 + c_0b_0|\lambda|^2\lambda + c_0^2|\lambda|^4 \end{array} \quad (19)$$

or

$$\begin{array}{l} (b_{-1} + b_0)\lambda^* + (c_{-1} + c_0)\lambda^{*2} \\ + (b_{-1} + b_0)\lambda + (b_{-1}^2 - b_0^2)|\lambda|^2 + (b_{-1}c_{-1} - b_0c_0)|\lambda|^2\lambda^* \\ + (c_{-1} + c_0)\lambda^2 + (c_{-1}b_{-1} - c_0b_0)|\lambda|^2\lambda + (c_{-1}^2 - c_0^2)|\lambda|^4 \\ < 0 \end{array} \quad (20)$$

Considering (7), (20) becomes

$$h(\lambda^* + \lambda) + \frac{h}{2}(c_{-1} - c_0)|\lambda|^2(\lambda^* + \lambda) < 0 \quad (21)$$

Equivalently

$$2\text{Real}(\lambda) + (c_{-1} - c_0)|\lambda|^2 \text{Real}(\lambda) < 0 \quad (22)$$

In (22), Real($\lambda$) < 0, $|\lambda|^2 > 0$. If $c_{-1} - c_0 > 0$, then (22) holds. Considering $c_0$ and $c_{-1}$ in (7), $c_{-1} - c_0 > 0$ is guaranteed by

$$\frac{1}{\omega_{select}^2} - \frac{h}{2\omega_{select}} \cot(\omega_{select} \frac{h}{2}) > 0 \quad (23)$$

or

$$\omega_{select} \frac{h}{2} \cot(\omega_{select} \frac{h}{2}) < 1 \quad (24)$$

The function $x\cot(x)$ is plotted as Fig. 5. It is learnt that if $0 < x < \pi$, then $x\cot(x) < 1$. Therefore if

$$0 < \omega_{select} \frac{h}{2} < \pi \quad (25)$$

then (24) holds. The step size is thus restrained by

$$0 < h < \frac{2\pi}{\omega_{select}} \quad (26)$$

to guarantee (24).

In summary, provided that the step size is smaller than one period defined by $\omega_{select}$, the numerical solution from Integrator A will follow the trend of the analytical solution despite the exact value of the step size. This is an interesting new type of numerical stability which is similar to being A-stable [27]-[30].

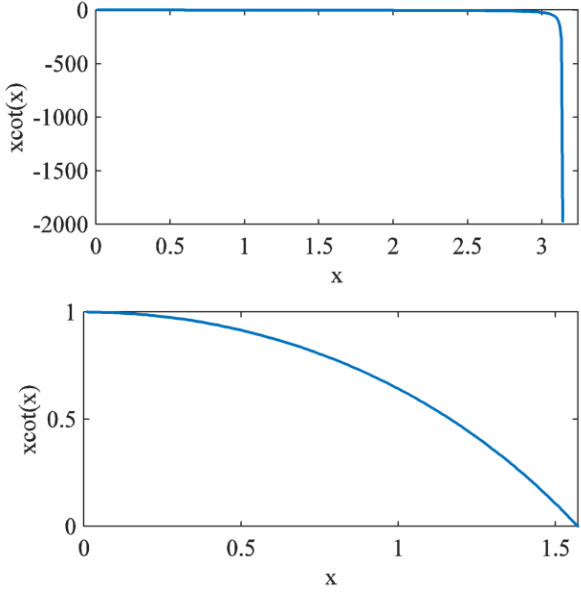

Fig. 5. Plot of *xcot*(*x*). Top: the plot drawn on the interval (0, π). Bottom: the plot drawn on the interval (0, $\frac{\pi}{2}$).

When Integrator A is applied, there is assumed to be a dominant frequency component other than 0 of interest in the system under study. Typically $\omega_{select}$ is chosen at the dominant angular frequency. Even for a very coarse depiction of a waveform dominated by this frequency, multiple samples in one period are required, which means a step size definitely smaller than the period. Therefore the constraint on the step size from the numerical stability aspect is in fact not a big issue.

*B. Numerical Stability of Integrator B*

Assume that the step size $h$ is smaller than half the period defined by $\omega_{select}$ so that $b_0 > 0$ and $c_0 < 0$ in (8). Applying Integrator B, the test ODE (11) is discretized as

$$x_t = x_{t-h} + b_0 \dot{x}_t + c_0 \ddot{x}_t = x_{t-h} + b_0 \lambda x_t + c_0 \lambda^2 x_t \quad (27)$$

For the numerical solution to tend to 0, it is required that

$$\left| \frac{1}{1 - b_0 \lambda - c_0 \lambda^2} \right| < 1 \quad (28)$$

which is guaranteed by

$$\left| 1 - b_0 \lambda - c_0 \lambda^2 \right| > 1 \quad (29)$$

guaranteed by

$$\left| 1 - b_0 \lambda - c_0 \lambda^2 \right|^2 > 1 \quad (30)$$

or

$$(1 - b_0 \lambda - c_0 \lambda^2)(1 - b_0 \lambda - c_0 \lambda^2)^* > 1 \quad (31)$$

namely

$$(1 - b_0 \lambda - c_0 \lambda^2)(1 - b_0 \lambda^* - c_0 \lambda^{*2}) > 1 \quad (32)$$

namely

$$(1 - b_0 \lambda^* - c_0 \lambda^{*2}) - b_0 \lambda (1 - b_0 \lambda^* - c_0 \lambda^{*2})$$
$$- c_0 \lambda^2 (1 - b_0 \lambda^* - c_0 \lambda^{*2}) \quad (33)$$
$$> 1$$

namely

$$-b_0 \lambda^* - b_0 \lambda + b_0^2 |\lambda|^2 + c_0^2 |\lambda|^4$$
$$-c_0 \lambda^2 - c_0 \lambda^{*2} + c_0 b_0 |\lambda|^2 \lambda + b_0 c_0 |\lambda|^2 \lambda^* \quad (34)$$
$$> 0$$

In (34),

$$-b_0 \lambda^* - b_0 \lambda + b_0^2 |\lambda|^2 + c_0^2 |\lambda|^4 > 0 \quad (35)$$

So it suffices to have the following to make (34) hold

$$-c_0 \lambda^2 - c_0 \lambda^{*2} + c_0 b_0 |\lambda|^2 \lambda + b_0 c_0 |\lambda|^2 \lambda^* > 0 \quad (36)$$

or

$$-\lambda^2 - \lambda^{*2} + b_0 |\lambda|^2 \lambda + b_0 |\lambda|^2 \lambda^* < 0 \quad (37)$$

namely

$$-2(\text{Real}^2(\lambda) - \text{Imag}^2(\lambda)) + 2b_0 |\lambda|^2 \text{Real}(\lambda) < 0 \quad (38)$$

In (38),

$$2b_0 |\lambda|^2 \text{Real}(\lambda) < 0 \quad (39)$$

So it suffices to have the following to make (38) hold

$$-2(\text{Real}^2(\lambda) - \text{Imag}^2(\lambda)) < 0 \quad (40)$$

Equivalently

$$|\text{Real}(\lambda)| > |\text{Imag}(\lambda)| \quad (41)$$

Note that (41) defines an infinite wedge on the left half-plane.

In summary, provided that the step size is smaller than half the period defined by $\omega_{select}$ and the eigenvalue $\lambda$ falls within the infinite wedge defined by (41), the numerical solution from Integrator B will follow the trend of the analytical solution despite the exact value of the step size. This is an interesting new type of numerical stability which is similar to being A(α)-stable or "nearly" A-stable [29]-[30].

Its numerical stability makes Integrator B unsuitable for being used as a main numerical integrator in time domain simulation. Some systems may contain eigenvalues which will cause a numerical blow-up in time domain simulation carried out with Integrator B. Nevertheless, its zero $b_{-1}$ and $c_{-1}$ coefficients make it suitable for dealing with discontinuities. As only a few time steps are to be calculated with Integrator B immediately after a discontinuity, if there is any, its numerical stability should not result in a numerical blow-up. In fact, this is how the novel power system transient simulation scheme [21] utilizes Integrator B.

## IV. TRANSIENT PERFORMANCE

In this section, a special case of (11) is considered where $\lambda$ is assumed to be real-valued with large magnitude. In this case, the analytical solution is a fast and monotonically decreasing signal tending to 0 [36]. The numerical solution from a numerical integrator is expected to follow the similar path to 0. However, some A-stable numerical integrators may create a numerical solution following a different path though finally tending to 0. For example, the implicit trapezoidal method induces the notorious numerical oscillation [32]-[33]. Such observation implies that A-stability itself is not enough to qualitatively characterize the accuracy of numerical integrators for fast and monotonically decreasing transients. The present paper thus defines the term "transient performance" to characterize the accuracy of numerical integrators from this aspect. It is desirable that a numerical integrator has such a property so that its numerical solution also monotonically and dramatically goes to 0 despite the step size.

Applying a numerical integrator considering second order derivative to the test ODE (11), the following recursive expression is obtained

$$x_t = \frac{1 + b_{-1}\lambda + c_{-1}\lambda^2}{1 - b_0\lambda - c_0\lambda^2} x_{t-h} \qquad (42)$$

The coefficient (43) determines the transient performance of the numerical integrator

$$\frac{1 + b_{-1}\lambda + c_{-1}\lambda^2}{1 - b_0\lambda - c_0\lambda^2} \qquad (43)$$

An A-stable numerical integrator merely assures that the magnitude of the coefficient (43) is smaller than 1 so that the numerical solution tends to 0. If the coefficient is close to -1, oscillating behavior is induced in the numerical solution. On the other hand, if the coefficient is close to 1, the numerical solution is monotonically decreasing but at a very slow rate. Neither situation works well for a rapidly decaying transient when the magnitude of $\lambda$ is large. Therefore numerical integrators exhibiting the above two types of transient performance cannot be understood as qualitatively accurate for fast and monotonically decreasing transients. Transient performance of the four numerical integrators studied in this paper is discussed as follows.

*A. Integrator A*

As discussed in Section III.A, when using Integrator A, it is assumed that the step size is smaller than one period defined by $\omega_{select}$. By (26) and Fig. 5, in (7)

$$\begin{aligned} c_0 &= -\frac{1}{\omega_{select}^2} + \frac{h}{2\omega_{select}} \cot(\frac{\omega_{select} h}{2}) \\ &= \frac{1}{\omega_{select}^2}(-1 + \frac{\omega_{select} h}{2} \cot(\frac{\omega_{select} h}{2})) \\ &< 0 \\ c_{-1} &= \frac{1}{\omega_{select}^2} - \frac{h}{2\omega_{select}} \cot(\frac{\omega_{select} h}{2}) \\ &= \frac{1}{\omega_{select}^2}(1 - \frac{\omega_{select} h}{2} \cot(\frac{\omega_{select} h}{2})) \\ &> 0 \end{aligned} \qquad (44)$$

Therefore the denominator in (43) must be larger than 1.

Consider

$$\frac{c_{-1}}{(\frac{h^2}{12})} = \frac{3}{(\frac{\omega_{select} h}{2})^2}(1 - \frac{\omega_{select} h}{2} \cot(\frac{\omega_{select} h}{2})) \qquad (45)$$

The plot of the function

$$\frac{3}{x^2}(1 - x\cot(x)) \qquad (46)$$

is given in Fig. 6, from which it is learnt that on the interval $(0, \pi)$, the function value is greater than 1, implying that

$$\frac{c_{-1}}{(\frac{h^2}{12})} = \frac{3}{(\frac{\omega_{select} h}{2})^2}(1 - \frac{\omega_{select} h}{2} \cot(\frac{\omega_{select} h}{2})) > 1 \qquad (47)$$

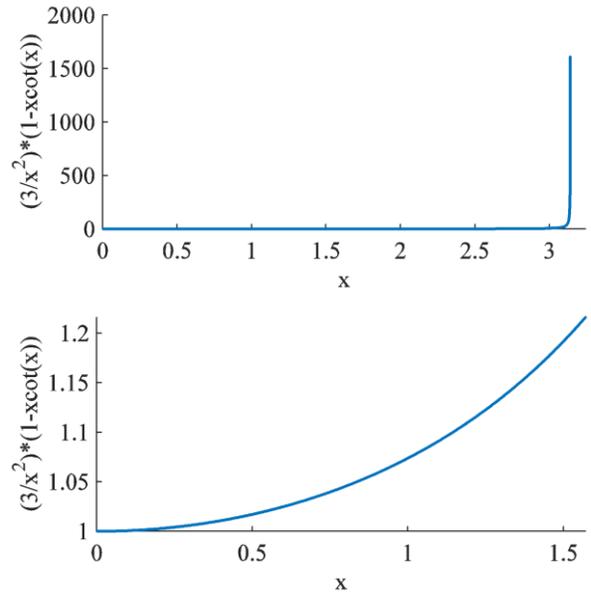

Fig. 6. Plot of $\frac{3}{x^2}(1-x\cot(x))$. Top: the plot drawn on the interval $(0, \pi)$. Bottom: the plot drawn on the interval $(0, \frac{\pi}{2})$.

$$c_{-1} > \frac{h^2}{12} \qquad (48)$$

Looking at the numerator of (43)

$$c_{-1}\lambda^2 + b_{-1}\lambda + 1 \qquad (49)$$

Note that

$$b_{-1}^2 - 4c_{-1} < (\frac{h}{2})^2 - 4\frac{h^2}{12} = -\frac{h^2}{12} < 0 \qquad (50)$$

The numerator thus must be larger than 0. Consequently (43) must be larger than 0 as well. No numerical oscillation will be induced by Integrator A. Nevertheless, $\lambda$ with a large magnitude makes (43) close to 1. The numerical solution is a slowly decaying transient, indicating the inaccuracy of Integrator A for rapidly decaying transients.

*B. Integrator B*

As discussed in Section III.B, when using Integrator B, it is assumed that the step size is smaller than half the period defined by $\omega_{select}$, namely

$$0 < h < \frac{\pi}{\omega_{select}} \qquad (51)$$

Considering (51) and (8), the denominator of (43) must be larger than 1; while the numerator is 1. Consequently (43) must be larger than 0. So no numerical oscillation will be induced by Integrator B. Moreover, $\lambda$ with a large magnitude makes (43) close to 0. The numerical solution is also a rapidly decaying transient, similar to the analytical solution. Therefore Integrator B for fast decaying transients is qualitatively accurate.

*C. Integrator C*

Substituting (9) into (43), the denominator must be larger than 1. The numerator (49) needs more discussion. Note that

$$b_{-1}^2 - 4c_{-1} = (\frac{h}{2})^2 - 4\frac{h^2}{12} = -\frac{h^2}{12} < 0 \quad (52)$$

The numerator must be larger than 0. Consequently (43) is larger than 0. No numerical oscillation will be induced by Integrator C. Unfortunately, $\lambda$ with a large magnitude makes (43) close to 1. As a result, the numerical solution is a slowly decaying transient, different from the analytical solution. Therefore Integrator C is inaccurate for fast decaying transients.

### D. Integrator D

Roughly following the same argument for Integrator B, Integrator D is shown to be qualitatively accurate for fast decaying transients in that the numerical solution follows the same path as the analytical solution tending to 0. Note that when using Integrator D, the step size is not confined.

## V. TIME DOMAIN VERIFICATION

### A. Numerical Error

Numerical error of Integrators A, B, C and D is studied in time domain by comparing the numerical solutions to the analytical solution. As a reference, the conventional implicit trapezoidal method (TR) and backward Euler method (BE) are also included.

Specifically, the following continuous linear system is to be used as the test system

$$\dot{x} = ax + bu \quad (53)$$

where $x$ is the state variable; $a = -5$; $b = 300$; the input $u$ is given as $u = cos(\omega_{syn}t)$, where $\omega_{syn} = 120\pi$ s$^{-1}$. Note that the derivative of the input is $\dot{u} = -\omega_{syn}sin(\omega_{syn}t)$. Further suppose that the initial value of $x$ is $x_0$. The analytical solution is [36]

$$x = (x_0 + \frac{ab}{\omega_{syn}^2 + a^2})e^{at} + b(\frac{\omega_{syn}}{\omega_{syn}^2 + a^2}sin(\omega_{syn}t)$$

$$- \frac{a}{\omega_{syn}^2 + a^2}cos(\omega_{syn}t)) \quad (54)$$

Numerical integrators are applied to the test system (53) so that the numerical solutions are obtained. The comparison between the numerical solutions and the analytical solution shows numerical error of the numerical integrators. For these numerical integrators, $\omega_{select}$ is set to $\omega_{syn}$ if applicable to make them accurate for the frequency component at $\omega_{syn}$.

In order to quantitatively study the numerical error, the following error measurement is used. The relative error regarding $x$ from a numerical integrator with a specified step size is defined as

$$err(x) = \frac{\|x_{num} - x_{ref}\|_2}{\|x_{ref}\|_2} \times 100 \quad (55)$$

where $x_{num}$ is the computed value by a numerical integrator with a specified step size; $x_{ref}$ is the reference value which comes from the analytical solution. As numerical solutions are discrete, the 2-norm is calculated at common time instants of $x_{num}$ and $x_{ref}$. Two case studies are presented as follows.

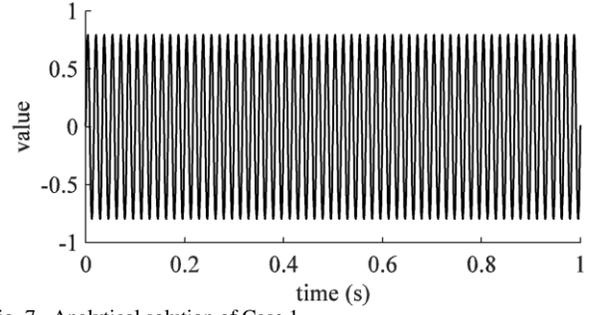

Fig. 7. Analytical solution of Case 1.

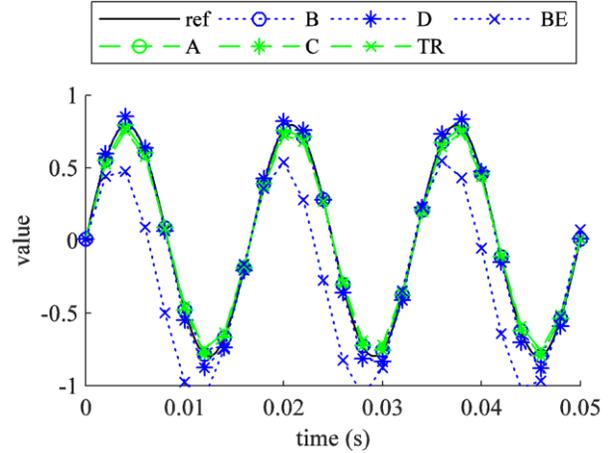

Fig. 8. Comparison of numerical solutions and analytical solution of Case 1.

TABLE I
ERROR OF NUMERICAL INTEGRATORS WITH DIFFERENT STEP SIZES IN CASE 1

| Step Size (μs) | Numerical Integrator | | | | | |
|---|---|---|---|---|---|---|
| | A | B | C | D | TR | BE |
| 125 | 0.0000 | 0.0000 | 0.0000 | 0.0370 | 0.0185 | 2.5803 |
| 250 | 0.0000 | 0.0000 | 0.0000 | 0.1480 | 0.0740 | 5.1598 |
| 500 | 0.0000 | 0.0000 | 0.0002 | 0.5920 | 0.2962 | 10.3179 |
| 1000 | 0.0000 | 0.0000 | 0.0028 | 2.3723 | 1.1870 | 20.6419 |
| 2000 | 0.0000 | 0.0000 | 0.0455 | 9.5852 | 4.7822 | 41.4123 |
| 4000 | 0.0000 | 0.0000 | 0.7593 | 40.1607 | 19.7071 | 84.2506 |

#### 1) Case 1: Sinusoidal Steady State

In Case 1, the initial value $x_0$ is set to

$$x_0 = -\frac{ab}{\omega_{syn}^2 + a^2} \quad (56)$$

so that the system directly enters sinusoidal steady state. Fig. 7 shows the analytical solution of this case. Numerical solutions from different numerical integrators with a step size of 2 ms are compared to the analytical solution in Fig. 8. Table I lists the error of the numerical integrators with different step sizes. Computations performed in Table I are from 0 to 1 s.

From Table I it is learnt that Integrators A and B are exactly accurate despite the step size, as is theoretically predicted in Section II. Integrator C is quite accurate; if the step size is doubled, the error roughly increases by a factor of 16. For Integrator D, if the step size is doubled, the error increases by a factor of around 4. The relation between the error and the step size of TR is the same as that of Integrator D, but it is about two times more accurate than Integrator D if the same step size is used. BE is the least accurate; if the step size is doubled, the

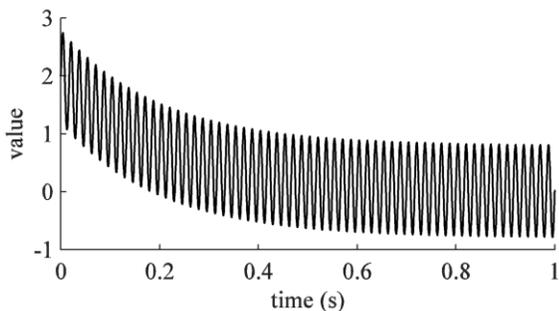
Fig. 9. Analytical solution of Case 2.

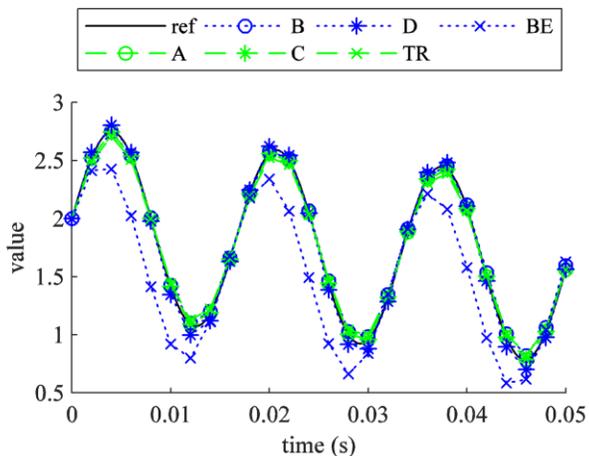
Fig. 10. Comparison of numerical solutions and analytical solution of Case 2.

TABLE II
ERROR OF NUMERICAL INTEGRATORS WITH DIFFERENT STEP SIZES IN CASE 2

| Step Size (μs) | Numerical Integrator | | | | | |
|---|---|---|---|---|---|---|
| | A | B | C | D | TR | BE |
| 125 | 0.0000 | 0.0194 | 0.0000 | 0.0245 | 0.0123 | 1.7052 |
| 250 | 0.0000 | 0.0774 | 0.0000 | 0.0980 | 0.0490 | 3.4093 |
| 500 | 0.0000 | 0.3100 | 0.0001 | 0.3921 | 0.1962 | 6.8152 |
| 1000 | 0.0000 | 1.2466 | 0.0019 | 1.5702 | 0.7857 | 13.6258 |
| 2000 | 0.0000 | 5.1240 | 0.0301 | 6.3369 | 3.1616 | 27.3049 |
| 4000 | 0.0001 | 23.2684 | 0.5010 | 26.4994 | 13.0036 | 55.4493 |

error is also doubled.

*2) Case 2: Time Domain Response with an Initial Transient*

In Case 2, $x_0$ is set to 2. As a result, there is an initial transient at the beginning. Fig. 9 shows the analytical solution of this case. Numerical solutions from different numerical integrators with a step size of 2 ms are compared to the analytical solution in Fig. 10. Table II lists the error of the numerical integrators with different step sizes. Computations performed in Table II are from 0 to 1 s.

From Table II it is learnt that Integrator A is highly accurate; it almost introduces no error. Integrator C is quite accurate introducing little error. For Integrator B, if the step size is doubled, the error increases by a factor of around 4. Integrator D is roughly as accurate as Integrator B. If the same step size is used, TR is about twice as accurate as Integrators B and D. Again, BE is the least accurate. The relation between the error and the step size of Integrators C and D, TR and BE in this case is the same as that in the previous case.

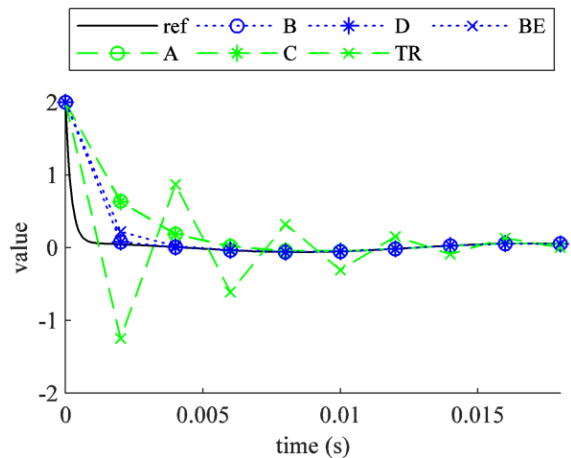
Fig. 11. Transient performance of numerical integrators.

*B. Transient Performance*

In order to better demonstrate the transient performance of the numerical integrators, the parameter *a* of the system (53) is changed to -5000; other conditions are the same as Case 2 in the previous subsection. Numerical solutions from the numerical integrators under consideration are compared to the analytical solution in Fig. 11 using a step size of 2 ms.

In Fig. 11, TR exhibits the notorious numerical oscillation. Transient performance of Integrators A, B, C and D is as theoretically predicted in Section IV. Transient performance of Integrators A and C is unsatisfactory as the damping effect is poor leading to significant mismatch between the numerical solutions and the analytical solution at the beginning. The damping effect of Integrators B and D is stronger than that of BE, resulting in a better match with the analytical solution.

One solution to the unsatisfactory transient performance of Integrators A and C is using Integrators B and D at the beginning for several time steps or half time steps to exploit their agreeable transient performance, and later switching back to Integrators A and C to exploit their high accuracy. This idea is motivated by the Critical Damping Adjustment (CDA) technique [37] broadly used in power system transient simulation [9], [38]-[39].

## VI. CONCLUSION AND FUTURE WORK

Comprehensive studies on frequency response optimized integrators considering second order derivative are presented in this paper, regarding numerical error, numerical stability and transient performance. Theoretical observations are verified by time domain case studies. Features of the numerical integrators under study are summarized as follows:

1) Integrator A is highly accurate for signals with a nonzero dominant frequency component as well as slow variants. Provided the step size is smaller than one period defined by $\omega_{select}$, Integrator A is qualitatively accurate for fast transients and high frequency oscillations despite the exact value of the step size. These signals or components can be skipped with relatively large step sizes in time domain simulation if not of interest. Therefore Integrator A is suitable to be used as a main numerical integrator for circuits or systems of which the signals inside are dominated by a specific frequency component, such as the carrier frequency in electronic communication circuits and the nominal fundamental frequency in power system networks. Unfortunately the ability of Integrator A to track fast decaying transients is weak. This drawback can be mitigated by temporarily using Integrator B immediately after discontinuities, which typically induce such transients.

2) Integrator B is accurate for signals with a nonzero dominant frequency component and slow variants. However some high frequency oscillations may cause numerical blow-up in time domain simulation carried out with Integrator B due to its unsatisfactory numerical stability. Consequently Integrator B is not a good choice to be used as a main numerical integrator in time domain simulation. Nevertheless, its strong ability to track fast decaying transients makes it a good complement to Integrator A for discontinuities.

3) Integrator C is highly accurate for signals around 0 Hz. Moreover it is A-stable. Its high accuracy and pleasant numerical stability make it suitable to be used as a main numerical integrator for the parts of the circuit or system under study where the signals inside are slow variants. However its ability to track fast decaying transients is weak so it has to be complemented by Integrator D immediately after discontinuities.

4) Integrator D is relatively accurate for slow variants. It has the pleasant numerical stability of being L-stable. In addition, its ability to track fast decaying transients is strong. Considering its advantages, it may be used as a main numerical integrator in time domain simulation. However, it is used as a complement to Integrator C instead in the novel power system transient simulation scheme [21] because the accuracy of Integrator C is much higher.

In the future, research effort may be directed to further development of frequency response optimized integrators theoretically and practically. On the theoretical side, relation between frequency domain error analysis and the conventional Taylor expansion-based truncation error analysis may be investigated; frequency response optimized integrators may be generalized to higher order derivatives and multiple time steps. On the practical side, applications to different areas and efficient implementation of these numerical integrators may be interesting topics.